\documentclass[12pt]{article}
\usepackage{latexsym}
\usepackage{graphicx}
\usepackage{cite}
\usepackage{enumitem}
\usepackage{amssymb,mathrsfs,amsmath}
\usepackage{booktabs}
\usepackage[justification=centering]{caption}
\usepackage{subfigure}
\usepackage{float}
\usepackage{color}
\usepackage{arydshln}

\usepackage{mathtools,amsthm}
\usepackage[normalem]{ulem}
\usepackage{mathabx}
\usepackage{acronym}

\usepackage[margin=1.5in]{geometry}%
\newtheorem{theorem}{Theorem}[part]

\newtheorem{lemma}{Lemma}[part]
\newtheorem{corollary}{Corollary}[part]

\newtheorem{example}{Example}[part]

\def \ep{\hbox{ }\hfill$\Box$}

\def\top{\mathsf{T}}

\begin{document}

\title{An existence result for weakly homogeneous variational inequalities}

\author{Meng-Meng Zheng\thanks{School of Mathematics, Tianjin University, Tianjin 300350, P.R. China.
Email: zmm941112@tju.edu.cn.}
\and Zheng-Hai Huang\thanks{Corresponding Author. School of Mathematics,
Tianjin University, Tianjin 300350, P.R. China. Email:
huangzhenghai@tju.edu.cn. Tel:+86-22-27403615 Fax:+86-22-27403615}}
\date{}%March 18, 2019}

\maketitle

\begin{abstract}
In this paper, what we concern about is the weakly homogeneous variational inequality over a finite dimensional real Hilbert space. We achieve an existence result {under} copositivity of leading term of the involved map, norm-coercivity of the natural map and several additional conditions. These conditions we used are easier to check and cross each other with those utilized in the main result established by Gowda and Sossa (Math Program 177:149-171, 2019). As a corollary, we obtain a result on the solvability of nonlinear equations with weakly homogeneous maps involved. Our result enriches the theory for weakly homogeneous variational inequalities and its subcategory problems in the sense that the main result established by Gowda and Sossa covers a majority of existence results on the subcategory problems of weakly homogeneous variational inequalities. Besides, we compare our {existence} result with the well-known coercivity result obtained for general variational inequalities and a norm-coercivity result obtained for general complementarity problems, respectively.
\end{abstract}

\noindent {\bf Key words:}\hspace{2mm} Weakly homogeneous map, variational inequality, copositivity, norm-coercivity.

\noindent {\bf Mathematics Subject
Classifications (2010):}\hspace{2mm} 65K10, 90C33 \vspace{3mm}

\section{Introduction}
Very recently, Gowda and Sossa investigated a class of continuous map, named as weakly homogeneous map, and the corresponding variational inequality {(VI)} over a finite dimensional real Hilbert space in \cite{GS-18}. As an application, they discussed the solvability of nonlinear equations {with} weakly homogeneous maps over closed convex cones, which covers tensor equations (or multilinear systems) \cite{Han-17,WCW-19,CLS-19,ZLC-20,HLQZ-18,XJW-18,DW-16} as special cases. Let $\mathbb{H}$ be a finite dimensional real Hilbert space with inner product $\langle \cdot,\cdot\rangle$ and norm $\|\cdot\|$, $C$ be a closed convex cone in $\mathbb{H}$ and $\mathbb{R}_{++}$ be the set of positively real numbers. Recall that a continuous map $h$: $C\rightarrow \mathbb{H}$ is called positively homogeneous of degree $\gamma(\geq0)$ if and only if $h(\lambda\mathbf{x})=\lambda^{\gamma}h(\mathbf{x})$ for any $\mathbf{x}\in C$ and $\lambda\in \mathbb{R}_{++}$. As a generalization of the positively homogeneous map, a weakly homogeneous map of degree $\gamma(\geq0)$ $f$ is defined to be a sum of a positively homogeneous of degree $\gamma(\geq0)$ $h$ and a remainder which is continuous on $C$ and satisfies $\lim\limits_{\|x\|\rightarrow\infty}\frac{f(x)-h(x)}{\|x\|^\gamma}=0$ where $x\in C$. Owing to the fact that $\lim\limits_{\lambda\rightarrow\infty}\frac{f(\lambda x)}{\lambda^\gamma}=h(x)$ for all $x\in C$, the positively homogeneous part $h$ is often called the ``leading term'' or ``recession map'' of the weakly homogeneous map $f$ and denoted by $f^\infty$.

Given a weakly homogeneous map $f$ on a closed convex cone $C$ and a closed convex subset $K$ in $C$, the  weakly homogeneous variational inequality (WHVI) {\cite{GS-18,MZH-20}}, denoted by $\mbox{\rm WHVI}(f,K)$ with $\mbox{\rm SOL}(f,K)$ being its solution set, is to find a vector
$x^*\in H$ such that
\begin{equation}\label{whvi}
x^*\in K,\quad \langle f(x^*),y-x^*\rangle\geq0,\quad \mbox{\rm for\ all}\ y\in K.
\end{equation}
When $K$ is a cone, \eqref{whvi} is called the weakly homogeneous complementarity problem (WHCP), denoted by $\mbox{\rm WHCP}(f,K)$.
WHVIs and WHCPs cover several recently researched special VIs and CPs as subcategory problems. Specifically, when $f$ is a polynomial, $\mbox{\rm WHVI}(f,K)$ and $\mbox{\rm WHCP}(f,K)$ come back to a polynomial variational inequality (PVI) studied in \cite{H-18} and polynomial complementarity problem (PCP) studied in \cite{G-17,LHL-18}, respectively; furthermore, when $f$ is a homogeneous polynomial, $\mbox{\rm WHVI}(f,K)$ and $\mbox{\rm WHCP}(f,K)$ come back to a tensor variational inequality (TVI) studied in \cite{WHQ-18} and tensor complementarity problem (TCP) studied in \cite{HQ-19-1,HQ-19-2,HQ-19-3}, respectively.

In \cite{GS-18}, the authors established some close connections between VIs and CPs with involved maps being weakly homogeneous map of positive degree. For instance, let $K^\infty:=\{u\in H:u+K\subseteq K\}$ represent the recession cone \cite{R70} of $K$, the main result given in \cite[Theorem 4.1]{GS-18} showed that if the corresponding recession cone CP $\mbox{\rm WHCP}(f^\infty, K^\infty)$ has and only has the only zero solution, and the (topological) index of the natural map \cite{FP03} of $\mbox{\rm WHCP}(f^\infty, K^\infty)$ at the origin is nonzero, then $\mbox{\rm WHVI}(f, K)$ has a nonempty, compact solution set. To our best knowledge, this degree-theoretic theorem is different from the famous coercivity result \cite{FP03} and covers a majority of existence results on the subcategory problems of WHVIs, including the well-known Karamardian's theorem \cite{K-76} for homogeneous maps on proper cones. However, an undesirable fact is that the degree-theoretic condition is often not easy to check.

Inspired by above, we aim to find {some easy-verified} conditions to guarantee the existence of solutions of WHVIs. Our paper is organized as follows. In Section 2, we show an alternative theorem for WHVIs. In Section 3, we establish nonemptiness and compactness of solution sets to WHVIs with the help of copositivity and norm-coercivity conditions. In Section 4, we compare our results with existing ones, including the main result given in \cite[Theorem 4.1]{GS-18} and the well-known coercivity result given for VIs. Finally, we sum up the conclusions in the last section.

\section{An alternative theorem}
In this section, we give an alternative theorem for WHVIs by making use of degree theory. Let $\Omega$ be a bounded open set in $\mathbb{H}$ {and} map $\phi:\mbox{\rm cl}{\Omega}\rightarrow \mathbb{H}$ be continuous where cl${\Omega}$ denote the closure of $\Omega$ and $p\in \mathbb{H}$. Recall that the topological degree of $\phi$ over $\Omega$ with respect to $p$ is well-defined if $p\notin \phi(\partial \Omega)$ where ${\partial\Omega}$ denote the boundary of $\Omega$, which is an integer used to judge the existence of a solution to the equation $\phi(x)=p$, denoted by $\mbox{\rm deg}(\phi,\Omega,p)$. The following properties of the topological degree plays important roles in our proof of alternative theorem for WHVIs.

%\begin{lemma}\label{lem2}{\rm (the excision property of the degree)}{\rm\cite{FP03}}
%Let $\Omega$ be a nonempty, bounded open subset in $H$, $\phi:\mbox{\rm cl}{\Omega}\rightarrow H$ be continuous and $p\notin \phi(\partial\Omega)$. Then $\mbox{\rm deg}(\phi,\Omega,p)=\mbox{\rm deg}(\phi,\Omega_1,p)$ for every open subset $\Omega_1$ of $\Omega$ such that $p\notin \phi(\Omega\setminus\Omega_1)$.
%\end{lemma}
\begin{lemma}\cite{FP03}\label{Theo-1}
Let $K$ be a closed convex set in $C$ and $f: C \rightarrow \mathbb{H}$ be continuous. If there exists a bounded open set $U$ with $\mbox{\rm cl}U\subseteq C$ such that $\mbox{\rm deg}(F^{nat}_{K},U,0)\neq 0$ where $F^{nat}_{K}(x):=x-\Pi_{K}(x-F(x))$ is the natural map of $\mbox{\rm VI}(f,K)$ with $F$ being a given continuous extension of $f$ and $\Pi_{K}(x)$ meaning the orthogonal projection of an $x\in H$ onto $K$, then $\mbox{\rm VI}(f,K)$ has a solution in $U$.
\end{lemma}

\begin{lemma}\label{lem1}{\cite{FP03}}
Let $\Omega$ be a nonempty, bounded open subset in $\mathbb{H}$. Then $\mbox{\rm deg}(\mathcal{H}(\cdot, t),\Omega,p(t))$ is independent of $t\in [0,1]$ for any two continuous maps $\mathcal{H} : \mbox{\rm cl} \Omega \times [0,1] \rightarrow \mathbb{R}^{n}$ and
$p : [0,1] \rightarrow \mathbb{R}^{n}$ such that
    $p(t)\notin \mathcal{H}(\partial\Omega, t)$ for any $t\in [0,1].$
\end{lemma}

By Lemma \ref{lem1}, we can see when the continuous map $\varphi:\mbox{\rm cl}{\Omega}\rightarrow \mathbb{H}$ satisfies $\varphi(x)=0$ if and only if $x=0$, then,
$\text{deg}(\varphi,\Omega',0)$ is invariant for any bounded open set $\Omega'$ containing 0 and contained in $\Omega$, and the common degree is written as $\text{ind}(\varphi,0)$.

In addition, we need the following result given in \cite{GS-18}

\begin{lemma}\label{lemma-1}\cite{GS-18}
Let $K$ be a closed convex set in cone $C$ and $f: C \rightarrow \mathbb{H}$ be a weakly homogeneous map of positive degree. If $f^\infty$ is copositive on $K^\infty$ and {$\mbox{\rm SOL}(f^{\infty},K^\infty)=\{0\}$},
then $\mbox{\rm WHVI}(f,K)$ has a nonempty, compact solution set and ${\rm deg}(F^{nat}_{K},{\Omega},0)\neq 0$ for any bounded open set ${\Omega}$ containing $\mbox{\rm SOL}(f,K)$.
\end{lemma}

Recall that a map $\psi:D\rightarrow \mathbb{H}$ is said to be copositive on $D\subseteq \mathbb{H}$ \cite{HP-90}, if $\langle\psi(x)-\psi(0), x\rangle\geq0$ holds for any $x\in D$. Now, we show an alternative theorem for WHVIs by employing copositivity of maps and above lemmas.
\begin{theorem}\label{lemma-2}
Let $K$ be a closed convex set in cone $C$ and $f: C \rightarrow \mathbb{H}$ be a weakly homogeneous map of positive degree. If $f^\infty$ is copositive on $K^\infty$, then either the $\mbox{\rm WHVI}(f,K)$ has a solution or there exists an unbounded sequence $\{x_k\}\subseteq K$ and a positive sequence $\{t_k\}\subseteq(0,1)$ such that for each $k$,
%Let $K$ be a closed convex set in cone $C$ and $f: C \rightarrow H$ and $g: H \rightarrow H$ be two weakly homogeneous maps. Suppose that  $g^{-1}(C)\subseteq C$, $g(x)-g^\infty(x)\in C$ as $\|x\|\rightarrow \infty$, $(g^{\infty})^{-1}(K^\infty)\subseteq C$ and $g^\infty$ is an injective map on $H$. If $f^\infty$ is copositive with respect to $g^\infty$ on $K^\infty$, then either the WHGVI$(f,g,K)$ has a solution or there exists an unbounded sequence $\{x_k\}$ and a positive sequence $\{t_k\}\subseteq(0,1)$ such that $g(x_k)\in K$ and for each $k$,
\begin{equation}\label{*}
\langle f^\infty(x_k)+t_k x_k+(1-t_k) (f(x_k)-f^\infty(x_k)), y-x_k\rangle\geq0,\quad \forall y\in K.
\end{equation}
\end{theorem}

\noindent{\bf Proof.} Let $F$ and $F^\infty$ be any given continuous extensions of $f$ and $f^\infty$, respectively. For the sake of contradiction, we assume that $\mbox{\rm SOL}(f,K)=\emptyset$ and $$\bigcup\limits_{0<t<1}\mbox{\rm SOL}(f^\infty+t\mathcal{I}+(1-t)(f-f^\infty),K)$$ is bounded, where $\mathcal{I}$ means the identity map from $\mathbb{H}$ into $\mathbb{H}$.
Then, consider the following homotopy map:
$$\mathcal{H}(x,t)=x-\Pi_K(x-(F^\infty(x)+tx+(1-t)(F(x)-F^\infty(x)))),\quad \forall\ (x,t)\in \mathbb{H}\times[0,1].$$
It is easy to see that $\mathcal{H}(\cdot,t)$ is just the natural map of WHVI$(f^\infty+t\mathcal{I}+(1-t)(f-f^\infty),K)$ for each $t\in[0,1]$. Denote the set of zeros of $\mathcal{H}(\cdot,t)$ by:
$$
\mathbb{Z}:=\{x\in \mathbb{H}\mid\mathcal{H}(x,t)=0\; \text{for some}\ t\in [0,1]\}.
$$
Since $\mbox{\rm SOL}(f,K)=\emptyset$, it follows that $\{x\in \mathbb{H}\mid\mathcal{H}(x,0)=0\}$ is bounded, which, together with another assumption, implies that $\{x\in \mathbb{H}\mid\mathcal{H}(x,t)=0\; \text{for some}\ t\in [0,1)\}$ is bounded. Now, we consider the set $\{x\in \mathbb{H}\mid\mathcal{H}(x,1)=0\}$. When $t=1$, $\mathcal{H}(x,t)$ becomes
$$\mathcal{H}(x,1)=x-\Pi_{K}(x-(F^\infty(x)+x)).$$
Since $f^\infty$ is copositive on $K^\infty$, it is not difficult to see that $\mbox{SOL}((f^\infty+\mathcal{I})^\infty,K^\infty)=\{0\}$ and $(f^\infty+\mathcal{I})^\infty$ is copositive on $K^\infty$. So from Lemma \ref{lemma-1}, it follows that $\mbox{\rm SOL}(f^\infty+\mathcal{I},K)$ is nonempty and compact and ${\rm deg}((F^\infty+\mathcal{I})^{nat}_{K},{\Omega},0)\neq 0$ for any bounded open set ${\Omega}$ containing $\mbox{\rm SOL}(f^\infty+\mathcal{I},K)$.  Hence, the set $\mathbb{Z}$ is uniformly bounded.

Let $\Omega\subseteq\mathbb{Z}$ be a bounded open set in $\mathbb{H}$, then for all $t\in [0,1]$,
$0\notin \mathcal{H}(\partial\Omega,t)$.
By Lemma \ref{lem1}, it follows that,
$$
\text{deg}(\mathcal{H}(\cdot,0),\Omega,0)=\text{deg}(\mathcal{H}(\cdot,1),\Omega,0)=\text{deg}((F^{\infty}+\mathcal{I})^{nat}_{K},\Omega,0)\neq 0,$$
which means that $\mbox{SOL}(f,g,K)$ is nonempty by Lemma \ref{Theo-1}. A contradiction!
\ep
%In \cite{L-13}, an existence result for GCPs was shown as below:
%\begin{lemma}\label{lemma}\cite{L-13}
%Let $K$ be a closed convex cone in $H$. Let $F$ be a continuous map from $K$ into $H$ and $g$ be an injective continuous map from $H$ into $H$. Then either the GCP$(F,g,K)$ has a solution or there exist an unbounded sequence ${x_k}$ and a positive sequence ${\tau_k}$ such that
%\begin{equation}\label{*}
%K\ni g(x_k)\perp F(x_k)+\tau_k g(x_k)\in K^*,\quad \forall\; k.
%\end{equation}
%\end{lemma}

\section{An existence result}
Recall that in \cite{GS-18}, the authors investigated the nonemptiness and compactness of the solution set of WHVI$(f,K)$ well, and established some good results based on a condition that SOL$(f^\infty,K^\infty)=\{0\}$, which is an important condition for deriving the boundedness of the solution set of WHVI$(f,K)$. Noting that the norm-coercivity \cite{FP03} of the natural map $F^{nat}_K$ is another useful condition to obtain the boundedness of the solution set of WHVI$(f,K)$. {Below}, we establish an existence result for WHVI$(f,K)$ with the help of the norm-coercivity of the natural map $F^{nat}_K$, rather than the assumption SOL$(f^\infty,K^\infty)=\{0\}$.
\begin{theorem}\label{Theorem-WHVI}
Let $K$ be a closed convex set in cone $C$ and $f: C \rightarrow \mathbb{H}$ be weakly homogeneous of positive degree with $0\in K$. Suppose that the following conditions hold:
%Let $K$ be a closed convex set in cone $C$ with $0\in K$, $f: C \rightarrow H$ and $g: H\rightarrow H$ be two weakly homogeneous maps of degrees $\gamma_1$ and $\gamma_2$, respectively. Suppose that $g$ satisfies that $g^{-1}(C)\subseteq C$, $g(x)-g^\infty(x)\in C$ as $\|x\|\rightarrow \infty$, $(g^{\infty})^{-1}(K^\infty)\subseteq C$ and $g^\infty$ is an injective map on $H$. If the following conditions hold:
\begin{itemize}
\item[$(\mbox{\rm a})$] $f^\infty$ is copositive on $K$;
\item[$(\mbox{\rm b})$] there exists no $c>0$ such that $-x=c(f(x)-f^\infty(x))$ {for any $x\in K$ satisfying $\|x\|\rightarrow\infty$};
 \item[$(\mbox{\rm c})$] $\lim\limits_{x\in K,\|x\|\rightarrow\infty}\|F^{nat}_{K}(x)\|=\infty$ and $\|f(x)-f^\infty(x)\|\leq\|F^{nat}_{K}(x)\|$ for any $x\in K$ satisfying $\|x\|\rightarrow\infty$, where $F^{nat}_{K}(x):=x-\Pi_{K}(x-F(x))$ is the natural map of $\mbox{\rm VI}(f,K)$ with $F$ being a given continuous extension of $f$;
\end{itemize}
then $\mbox{\rm WHVI}(f,K)$ has a nonempty, compact solution set.
%then WHGVI$(f+p,g,K)$ has a nonempty, compact solution set for all $p\in H$.
\end{theorem}

\noindent{\bf Proof.} First, we show that $f^\infty$ is copositive on $K^\infty$. For any $u\in K^\infty$, there exist two sequences $\{t_k\}\subseteq \mathbb{R}_+$ and $\{x_k\}\subseteq K$ such that $u=\lim_{k\rightarrow\infty}\frac{x_k}{t_k}$. Thus, we have that for any $u\in K^\infty$,
$$\langle f^\infty(u), u\rangle=\lim_{k\rightarrow\infty}\langle f^\infty(\frac{x_k}{t_k}), \frac{x_k}{t_k}\rangle=\lim_{k\rightarrow\infty}\frac{\langle f^\infty(x_k), x_k\rangle}{t_k^{\gamma+1}}.$$
Noting that $f^\infty$ is copositive on $K$, which means that $\langle f^\infty(x_k), x_k\rangle\geq 0$, thus we have that $\langle f^\infty(u), u\rangle\geq 0$ for any $u\in K^\infty$. That is to say, $f^\infty$ is copositive on $K^\infty$.

{Second, we show that $\mbox{\rm SOL}(f,K)\neq\emptyset$. Suppose that $\mbox{\rm SOL}(f,K)=\emptyset$, then} from Theorem \ref{lemma-2}, there exists an unbounded sequence $\{x_k\}$ and a positive sequence $\{t_k\}\subseteq(0,1)$ such that $x_k\in K$ and \eqref{*} holds for each $k$. Noting that $0\in K$, thus we have that
$$
\langle x_k, f^\infty(x_k)+t_k x_k+(1-t_k) (f(x_k)-f^\infty(x_k))\rangle\leq0,
$$
which, together with the condition that $f^\infty$ is copositive on $K$ and $x_k\in K$, implies that
$$\langle x_k, t_k x_k+(1-t_k) (f(x_k)-f^\infty(x_k))\rangle\leq0.$$
Thereby, we can obtain that
\begin{equation}\label{equa-1}
t_k\|x_k\|\leq (1-t_k) \|f(x_k)-f^\infty(x_k)\|,
\end{equation}
{which implies that $f(x_k)-f^\infty(x_k)\neq0$ for sufficiently large $k$.
Furthermore, by $t_k>0$, condition (b) and the trigonometric inequality of the norm, we have that
\begin{equation}\label{123}
\|-t_k x_k+t_k(f(x_k)-f^\infty(x_k))\|<t_k\|x_k\|+t_k\|f(x_k)-f^\infty(x_k)\|
\end{equation}
for sufficiently large $k$.} In addition, by Lemma \ref{Theo-1} and \eqref{*}, it follows that
$$x_k=\Pi_K(x_k-f^\infty(x_k)-t_k x_k-(1-t_k) (f(x_k)-f^\infty(x_k)),$$
which, together with $F(x_k)=f(x_k)$ since $x_k\in K$, implies that {for sufficiently large $k$,}
\begin{equation}\label{111}
\begin{array}{lcl}
\|F^{nat}_{K}(x_k)\|
&=&\|x_k-\Pi_K(x_k-F(x_k))\|\\
&=&\|\Pi_K(x_k-f^\infty(x_k)-t_k x_k-(1-t_k) (f(x_k)-f^\infty(x_k)))-\Pi_K(x_k-f(x_k))\|\\
&\leq &\|-t_k x_k+t_k (f(x_k)-f^\infty(x_k))\|\\
&<& t_k\|x_k\|+t_k\|f(x_k)-f^\infty(x_k)\|\\
&\leq &(1-t_k) \|f(x_k)-f^\infty(x_k)\|+t_k\|f(x_k)-f^\infty(x_k)\|\\
&=&\|f(x_k)-f^\infty(x_k)\|,
\end{array}
\end{equation}
where the first inequality follows from non-expansiveness of Euclidean projector, the second inequality follows from \eqref{123} and the third inequality follows from \eqref{equa-1}.

On one hand, it follows from \eqref{111} that
$$\|F^{nat}_{K}(x_k)\|<\|f(x_k)-f^\infty(x_k)\|$$
{for sufficiently large $k$}. On the other hand, by condition (c), $x_k\in K$ and $\|x_k\|\rightarrow\infty$ as $k\rightarrow\infty$, we have $$\|f(x_k)-f^\infty(x_k)\|\leq\|F^{nat}_{K}(x_k)\|$$
{for sufficiently large $k$.} A contradiction yields! Hence, it follows that SOL$(f,K)\neq\emptyset$.

{Last}, from $\lim\limits_{x\in K,\|x\|\rightarrow\infty}\|F^{nat}_{K}(x)\|=\infty$, it follow that SOL$(f,K)$ is bounded. Therefore, we obtain that $\mbox{\rm SOL}(f,K)$ is nonempty and compact.
\ep

{In \cite{GS-18}, an extension of the $Z$-property for nonlinear maps was given as follows:
\begin{lemma}\label{Z}
Suppose $C$ is a closed convex cone with the dual cone being $C^*:=\{y\in \mathbb{H}\mid \langle y,x\rangle\geq0,\forall x\in C\}$ and $f: C\rightarrow\mathbb{H}$ satisfies:
\begin{equation}\label{000}
x\in C,y\in C^*\;\mbox{\rm and}\;\langle x,y\rangle=0\quad\Longrightarrow\quad\langle f(x),y\rangle\leq0.
\end{equation}
Then, $f(x^*)=q$ if and only if $q\in C$ and $x^*\in \mbox{\rm SOL}(\bar{f},C)$
where $\bar{f}(x):=f(x)-q$.
\end{lemma}

With Lemma \ref{Z}, we can obtain the following corollary of Theorem \ref{Theorem-WHVI}.
\begin{corollary}
Suppose conditions of Theorem \ref{Theorem-WHVI} hold with $K=K^\infty=C$ and $f$ satisfies \eqref{000}. Then for all $q\in C$, the equation $f(x)=q$ has a solution in $C$.
\end{corollary}}

\section{Some comparisons between our existence result and three related results}

Recall that for CPs, a related existence result which requires the norm-coercivity
of the natural map is as follows:

\begin{theorem}\label{Theorem-CP}\cite{FP03}
Let $K$ be a closed convex cone in $\mathbb{R}^n$ and $f: K\rightarrow\mathbb{R}^n$ be a continuous map. Suppose that for any $x\in K$, $\lim\limits_{\|x\|\rightarrow\infty}\|f_K^{nat}(x)\|=\infty$ and $\langle x, f(x)-f(0)\rangle\geq 0$, then $\mbox{\rm CP}(f,K)$ has a nonempty, compact solution set.
\end{theorem}

Actually, the condition that ``$\langle x, f(x)-f(0)\rangle\geq 0$'' implies that $f$ is a $q$-copositive map on cone $K$ with $q=f(0)$. From \cite[Theorem 5.1]{MZH-20}, it follows that if $f$ is a $q$-copositive map on cone $K$, then the condition (a) of Theorem \ref{Theorem-WHVI} must hold. Especially, when $K$ is a closed convex cone and the involved map $f(x)=f^\infty(x)+p$ with $p$ being a vector in $\mathbb{H}$, condition (a) in Theorem \ref{Theorem-WHVI} is equivalent to the statement that ``$\langle x, f(x)-f(0)\rangle\geq 0$'',  condition (b) in Theorem \ref{Theorem-WHVI} naturally holds and condition (c) in Theorem \ref{Theorem-WHVI} reduces to ``$\lim\limits_{x\in K,\|x\|\rightarrow\infty}\|f_K^{nat}(x)\|=\infty$''. That is to say, Theorem \ref{Theorem-WHVI} and Theorem \ref{Theorem-CP} coincide in the aforementioned case. However, for the case where $f-f^\infty$ is not a constant vector in $\mathbb{H}$, Theorem \ref{Theorem-WHVI} is different from Theorem \ref{Theorem-CP}, which can be illustrated from the following example, where WHCP$(f,K)$ satisfies all the conditions in Theorem \ref{Theorem-WHVI}, but don't satisfy the conditions of Theorem \ref{Theorem-CP}.

\begin{example}\label{Example-3.3}
{Consider $\mbox{\rm WHCP}(f,K)$ where $K=\{x\in\mathbb{R}^2: x_1\geq0,x_2\leq0\},$ and for any $x\in \mathbb{R}^2$,
$$f(x)=(x_1^2+x_2^2-\sqrt[4]{x_1^2}{+1},
-(x_1^2+x_2^2)+\sqrt[4]{x_2^2}{+2})^\top.$$}
Obviously, $K$ is a convex cone, and $f$ is weakly homogeneous of degree 2.
\begin{itemize}
\item For any $x\in K$, we have that
$\langle f^\infty(x), x\rangle=(x_1^2+x_2^2)(x_1-x_2)\geq0,$
thus $f^\infty$ is copositive on $K$.
\item For any $x\in K$ {with $\|x\|\rightarrow\infty$}, it is easy to see that there exists no $c>0$ such that $$-K\ni-x=c(f(x)-f^\infty(x))=c(-\sqrt{x_1}{+1},\sqrt{-x_2}{+2})^\top\in K.$$
\item For any $x\in K$, as $\|x\|\rightarrow\infty$, it follows that
    $$\left\{\begin{array}{l}
    (x-f(x))_1=x_1-(x_1^2+x_2^2)+\sqrt{x_1}{-1}\rightarrow-\infty\\ (x-f(x))_2=x_2+(x_1^2+x_2^2)-\sqrt{-x_2}{-2}\rightarrow+\infty.
    \end{array}\right.$$
So,
    %$x\in K^\infty\subseteq K$ and
%    \begin{itemize}
%    \item if $x_1x_2=0$, then
%    $x-f(x)\rightarrow(-\infty,0)^\top\in K^\infty$ or $x-f(x)\rightarrow(0,+\infty)^\top\in K^\infty$;
%    \item else $x-f(x)\rightarrow (+\infty,-\infty)^\top$.
%    \end{itemize}
%
    %$$\frac{\|f^{nat}_{K}(x)\|}{\|x\|}=\frac{\|x-\Pi_K(x-f(x))\|}{\|x\|}=\frac{\|f(x)\|}{\|x\|}>0.$$
    %or
$\Pi_K(x-f(x))=0$, $\lim\limits_{x\in K,\|x\|\rightarrow\infty}\|f^{nat}_{K}(x)\|=\lim\limits_{x\in K,\|x\|\rightarrow\infty}\|x-\Pi_K(x-f(x))\|=\infty$ and $$\lim\limits_{x\in K,\|x\|\rightarrow\infty}\frac{\|f^{nat}_{K}(x)\|}{\|f(x)-f^\infty(x)\|}=\lim_{x\in K,\|x\|\rightarrow\infty}\frac{\sqrt{x_1^2+x_2^2}}{\sqrt{(1-\sqrt{x_1})^2+(\sqrt{-x_2}+2)^2}}=\infty.$$
\end{itemize}
Thus, all conditions of Theorem \ref{Theorem-WHVI} are satisfied, which means that $\mbox{\rm SOL}(f,K)$ is nonempty and compact.

%Noting that
%$\text{SOL}(f^{\infty},K^\infty)=\{x\in K:x_1+x_2=0\}\neq\{0\},$
%thus the conditions of the main result given in \cite{GS-18} don't hold. Besides,
By taking $x\in K$ with $x_1=\frac{1}{4}$ and $x_2=0$, we have that
$\langle x, f(x)-f(0)\rangle=\frac{1}{64}-\frac{1}{8}<0,$
so the conditions in Theorem \ref{Theorem-CP} don't hold.
\end{example}

%Besides, recall that very recently in \cite{GS-18}, the authors achieved many good theoretical results on the nonemptiness and compactness of the solution set to WHVI$(f,K)$. In special, they gave a very wide existence result for WHVIs as the main result shown in \cite[Theorem 4.1]{GS-18} under the conditions that ind$\left((F^{\infty})^{nat}_{K^\infty}(x),0\right)\neq 0$ with $F^{\infty}$ being a continuous extension of $f^\infty$ and SOL$(f^{\infty},K^\infty)=\{0\}$, which covers a majority of existence results on the subcategory problems of WHVIs, including the Karamardian-type one in \cite[Theorem 5.1]{GS-18} and the copositivity one in \cite[Theorem 6.1]{GS-18}.
%Since Corollary \ref{Coro-WHCP} requires $0\in K$, which is unnecessary for \cite[Theorem 4.1]{GS-18}, thus the conditions in Corollary \ref{Coro-WHCP} cannot cover the conditions of \cite[Theorem 4.1]{GS-18}.

Second, we compare our result with the main result given in {\cite{GS-18}}:
\begin{theorem}\label{mainresult}\cite{GS-18}
Let $K$ be a closed convex set in cone $C$ and $f: C \rightarrow \mathbb{H}$ be a weakly homogeneous map of positive degree. If $\mbox{\rm SOL}(f^{\infty},K^\infty)=\{0\}$ and $\mbox{\rm ind}(G_{K^\infty},0)\neq 0$ where $G_{K^\infty}(x):=x-\Pi_{K^\infty}(x-G(x))$ with $G$ being a given continuous extension of $f^\infty$, then $\mbox{\rm WHVI}(f,K)$ has a nonempty, compact solution set.
\end{theorem}

Below, we construct three examples to show that conditions of Theorem \ref{Theorem-WHVI} and {Theorem \ref{mainresult}} cross each other.

\begin{example}\label{Example-3.4}
Consider $\mbox{\rm WHVI}(f,K)$ where $K:=\{x\in\mathbb{R}^2: x_1=x_2\geq-1\}$ and for any $x\in \mathbb{R}^2$, $$f(x)=((x_1-x_2)^2+x_1,(x_1-x_2)^2+x_2)^\top.$$
Obviously, $K$ is a convex subset in $\mathbb{R}^2$ with $0\in K$, and $f$ is weakly homogeneous of degree 2.
\begin{itemize}
\item For any $x\in K$, we have that
$\langle f^\infty(x), x\rangle=(x_1-x_2)^2(x_1+x_2)=0,$
thus $f^\infty$ is copositive on $K$.
\item For any $x\in K$ {with $\|x\|\rightarrow\infty$}, there exists no $c>0$ such that $$-x=c(f(x)-f^\infty(x))=cx.$$
\item For any $x\in K$, it follows that $x-f(x)=0$.
     Thus, we have that for any $x\in K$, $$\|f(x)-f^\infty(x)\|=\|f^{nat}_{K}(x)\|=\|x\|\quad \mbox{\rm and}\quad \lim\limits_{x\in K,\|x\|\rightarrow\infty}\|f^{nat}_{K}(x)\|=\lim\limits_{x\in K,\|x\|\rightarrow\infty}\|x\|=\infty.$$
\end{itemize}
Thus, all conditions of Theorem \ref{Theorem-WHVI} are satisfied, which means that $\mbox{\rm SOL}(f,K)$ is nonempty and compact.

However, noting that
$$\mbox{\rm SOL}(f^{\infty},K^\infty)=K^\infty=\{x\in\mathbb{R}^2: x_1=x_2\geq0\}\neq\{0\},$$
thus at least there is one condition of the main result given in {Theorem \ref{mainresult}} does not hold. That is to say, conditions in Theorem \ref{Theorem-WHVI} cannot be covered by those of {Theorem \ref{mainresult}}.
{In addition, we can see that conditions in Theorem \ref{Theorem-WHVI} are easy to check, but  it is not clear whether ${\rm deg}(G_{K^\infty},{\Omega},0)\neq 0$ holds or not where ${\Omega}\supseteq\mbox{\rm SOL}(f^\infty,K^\infty)$ is a bounded open set and $G_{K^\infty}(x)=x-\Pi_{K^\infty}(x-G(x))$ with $G$ being a given continuous extension of $f^\infty(x)=((x_1-x_2)^2,(x_1-x_2)^2)^\top$.}
\end{example}

\begin{example}\label{Example-3.5}
{Consider $\mbox{\rm WHCP}(f,K)$ where $K:=\mathbb{R}^2_+$ and for any $x\in \mathbb{R}^2$,
$$f(x)=((x_1+x_2)^2+{x_1},(x_1+x_2)^2+{x_2})^\top.$$}
Obviously, $K$ is a convex cone, and $f$ is weakly homogeneous of degree 2.
\begin{itemize}
\item For any $x\in K$, we have
$
\langle f^\infty(x), x\rangle=(x_1+x_2)^3\geq0,
$
thus $f^\infty$ is copositive on $K$.
\item For any $x\in K\setminus\{0\}$, there exists no $c>0$ such that $-x=c(f(x)-f^\infty(x))=cx$.
\item For any $x\in K$, it follows that
    $
    (x-f(x))_1=-(x_1+x_2)^2\rightarrow-\infty$ and $(x-f(x))_2=-(x_1+x_2)^2\rightarrow-\infty$ as $\|x\|\rightarrow\infty$,
    which implies that $\Pi_K(x-f(x))=0$ as $\|x\|\rightarrow\infty$.
    %$x\in K^\infty\subseteq K$ and
%    \begin{itemize}
%    \item if $x_1x_2=0$, then
%    $x-f(x)\rightarrow(-\infty,0)^\top\in K^\infty$ or $x-f(x)\rightarrow(0,+\infty)^\top\in K^\infty$;
%    \item else $x-f(x)\rightarrow (+\infty,-\infty)^\top$.
%    \end{itemize}
%
    %$$\frac{\|f^{nat}_{K}(x)\|}{\|x\|}=\frac{\|x-\Pi_K(x-f(x))\|}{\|x\|}=\frac{\|f(x)\|}{\|x\|}>0.$$
    %or
So, for any $x\in K$ and $\|x\|\rightarrow\infty$,
$$\|f(x)-f^\infty(x)\|=\|f^{nat}_{K}(x)\|=\|x\|\quad\mbox{\rm and}\quad \lim\limits_{x\in K,\|x\|\rightarrow\infty}\|f^{nat}_{K}(x)\|=\lim\limits_{x\in K,\|x\|\rightarrow\infty}\|x\|=\infty.$$
\end{itemize}
Thus, all conditions of Theorem \ref{Theorem-WHVI} are satisfied.

In addition,
%from \eqref{3.5},
it is not difficult to see that $\mbox{\rm SOL}(f^{\infty},K^\infty)=\{0\},$ which, together with $f^\infty$ is copositive on $K^\infty$, implies that conditions in {Theorem \ref{mainresult}} hold by Lemma \ref{lemma-1}.
\end{example}

\begin{example}
Consider $\mbox{\rm WHVI}(f,K)$ where $K=\{x\in\mathbb{R}^2: x_1\geq0,x_2\geq0,x_1+x_2\geq1\},$ and for any $x\in \mathbb{R}^2$,
$$f(x)=((x_1+x_2)^2-\sqrt[4]{x_1^2},
(x_1+x_2)^2-\sqrt[4]{x_2^2})^\top.$$

Obviously, $K$ is a convex subset in $\mathbb{R}^2$ with $K^\infty=\mathbb{R}^2_+$, and $f$ is weakly homogeneous of degree 2.

By a similar discussion as Example \ref{Example-3.5}, we can see that all conditions of {Theorem \ref{mainresult}} are satisfied. However, noting that $0\notin K$, thus those conditions in Theorem \ref{Theorem-WHVI} do not hold.
\end{example}

{It was shown in \cite{HP-90} that for a general VI$(f,K)$, employing the coercivity property of $f$ to establish existence is a common approach.} Last, we compare our result with the following well-known coercivity result for VIs.
\begin{theorem}\label{coercivity}{\cite{FP03}}
Let $K$ be a closed convex set in $\mathbb{R}^n$ and $f: K \rightarrow \mathbb{R}^n$ be continuous. If $f$ is coercive on $K$, that is there exists some $x^{ref}\in K$, $c>0$ and $\xi\geq0$ such that
$\langle f(x),x-x^{ref}\rangle\geq c\|x\|^\xi$ for any $x\in K$ with $\|x\|\rightarrow\infty,$
then VI$(f,K)$ has a nonempty, compact solution set.
\end{theorem}

The following example illustrates that the conditions in Theorem \ref{Theorem-WHVI} cannot be deduced by the above coercive one.

\begin{example}
{Consider $\mbox{\rm WHCP}(f,K)$ where $K:=\{x\in \mathbb{R}^2: x_1=x_2\}$, and for any $x\in \mathbb{R}^2$,
$$f(x)=(x_2^3-2x_1,-x_1^3)^\top.$$}
Obviously, $K$ is a convex cone, and $f$ is weakly homogeneous of degree 3.

First, we show that all conditions of Theorem \ref{Theorem-WHVI} are satisfied.
\begin{itemize}
\item For any $x\in K$, i.e., $x_1=x_2$ we have that
$\langle f^\infty(x), x\rangle=x_1(x_2^3)+x_2(-x_1^3)=0,$
thus $f^\infty$ is copositive on $K$.
\item For any $x\in K\setminus\{0\}$, there exists no $c>0$ such that $$(-x_1,-x_1)^\top=-x=c(f(x)-f^\infty(x))=(-2cx_1,0).$$
\item For any $x\in K$, it follows that $x-f(x)=(-x_1^3+3x_1,x_1^3+x_1)^\top$.
     Thus, we have that $\Pi_{K}(x-f(x))=(2x_1, 2x_1)$, for any $x\in K$, $$\|f(x)-f^\infty(x)\|=\|f^{nat}_{K}(x)\|=\|x\|\quad\mbox{\rm and}\quad\lim\limits_{x\in K,\|x\|\rightarrow\infty}\|f^{nat}_{K}(x)\|=\lim\limits_{x\in K,\|x\|\rightarrow\infty}\|-x\|=\infty.$$
\end{itemize}
Thus, all conditions of Theorem \ref{Theorem-WHVI} are satisfied, which means that $\mbox{\rm SOL}(f,K)$ is nonempty and compact.

Second, we show the coercivity condition of Theorem \ref{coercivity} does not hold. For any $x^{ref}\in K$, we have that $x_1^{ref}=x_2^{ref}$. Then for any $x\in K$, i.e., $x_1=x_2$, we have that
$$\langle f(x), x-x^{ref}\rangle=-2x_1^2-2x_1x_1^{ref}\rightarrow-\infty\quad\mbox{\rm as}\quad\|x\|\rightarrow\infty.$$
Thus, we cannot find some $x^{ref}\in K$, $c>0$ and $\xi\geq0$ such that
     $\langle f(x),x-x^{ref}\rangle\geq c\|x\|^\xi$ for any
   $x\in K$ satisfying $\|x\|\rightarrow\infty,$
which implies that $f$ is not coercive on $K$.
\end{example}

\section{Conclusions}
In this paper, we established an existence result for WHVI$(f,K)$ under the copositivity of $f^\infty$ on $K$, the norm-coercivity of the natural map $F^{nat}_K$ and some additional easy-verified conditions. Some examples were constructed to demonstrate that our conditions cannot deduced by the existing ones, especially the wide degree-theoretic theorem given for WHVIs in the main result of \cite{GS-18} and the well-known coercivity result given for general VIs in \cite{FP03}. Besides, our result also provided a supplement to the existence theory for those subclasses of WHVIs, such as PVIs and PCPs.

\end{document}